\documentclass[12pt]{article}
\begin{document}
\newcommand{\e}{{\bf{E}}}
\newcommand{\bea}{\begin{eqnarray}}
\newcommand{\ena}{\end{eqnarray}}
\newcommand{\beas}{\begin{eqnarray*}}
\newcommand{\enas}{\end{eqnarray*}}
\newcommand{\beq}{\begin{equation}}
\newcommand{\enq}{\end{equation}}
\newcommand{\G}{F}
\newcommand{\blambda}{\mbox{\boldmath {$\lambda$}}}
\def\qed{\hfill \mbox{\rule{0.5em}{0.5em}}}
\newtheorem{theorem}{Theorem}[section]
\newtheorem{example}{Example}[section]
\newtheorem{corollary}{Corollary}[section]
\newtheorem{proposition}{Proposition}[section]
\newtheorem{definition}{Definition}[section]
\newtheorem{lemma}{Lemma}[section]
\newtheorem{condition}{Condition}[section]

\begin{center}
{\Large \bf Distributional transformations, orthogonal polynomials, and Stein characterizations}

\vspace{7mm} Larry Goldstein\footnote{Department of Mathematics,
University of Southern California, Los Angeles CA 90089-2532,
USA}$^3$ and Gesine Reinert\footnote{Department of Statistics,
University of Oxford, Oxford OX1 3TG, UK, supported in part by
EPSRC grant no. GR/R52183/01}\footnote{This work was partially
completed while the authors were visiting the Institute for
Mathematical Sciences, National University of Singapore, in 2003.
The visit was supported by the Institute.}

\end{center}

\begin{abstract}
A new class of distributional transformations is introduced,
characterized by equations relating function weighted expectations
of test functions on a given distribution to expectations of the
transformed distribution on the test function's higher order
derivatives. The class includes the size and zero bias
transformations, and when specializing to weighting by polynomial
functions, relates distributional families closed under
independent addition, and in particular the infinitely divisible
distributions, to the family of transformations induced by their
associated orthogonal polynomial systems. For these families,
generalizing a well known property of size biasing, sums of
independent variables are transformed by replacing summands chosen
according to a multivariate distribution on its index set by
independent variables whose distributions are transformed by
members of that same family. A variety of the transformations
associated with the classical orthogonal polynomial systems have
as fixed points the original distribution, or a member of the same
family with different parameter.
\end{abstract}

\section{Introduction}

The zero bias transformation was introduced in \cite{goldrein}.
This mapping enjoys properties similar to those of the well known
size biased transformation (see e.g. \cite{gr96}) on non-negative
variables, but can be applied to mean zero random variables. One
main feature of the zero bias transformation is that its unique
fixed point is the mean zero normal distribution, and for this
reason it has been applied in Stein's method for the purpose of
normal approximation (\cite{goldrein}, \cite{hie}, \cite{zsm}, and
\cite{goldrein2}). The zero bias transformation is also related to
the $K_i$ function in the work \cite{hochen}, \cite{boldi},
\cite{st86}, and \cite{cacetal} and others; for a good overview
see \cite{chenshao}.

We place the classical size bias transformation and the zero bias
transformation in a broader context, showing that both are a
particular case of transforming a given distribution $X$ into
$X^{(P)}$ through the use of a measurable `biasing' function $P$.
To be more precise, for the given $X$ and $P$ let ${\cal C}^m$
denote the collection of functions whose $m^{th}$ derivative
exists and is measurable on ${\bf R}$, and suppressing $X$ on the
left hand side, set
$$
{\cal F}^m(P)=\{F \in {\cal C}^m: E|P(X)F(X)| < \infty \}.
$$
We consider transformations characterized by
\bea
\label{Gbias-intro} E P(X) F(X) = \alpha E F^{(m)}(X^{(P)}) \quad
\mbox{for all $F \in {\cal F}^m(P)$,}
\ena
where necessarily $\alpha = (m!)^{-1}EP(X) X^m$ when $X^m \in
{\cal F}(P)$; we insist $\alpha>0$. We coin this distribution the
$X-P$ biased distribution. For discrete distributions the
differential operator is replaced by the difference operator; see
Sections \ref{charlie} and \ref{kraw}.

Theorem \ref{biastheorem} provides our most general conditions on
the existence of transformations characterized by
(\ref{Gbias-intro}), where the `biasing' function $P$ is only
required to have $m$ sign changes and satisfy certain
orthogonality and positivity conditions; we call $m$ the order of
the resulting transformation. In the particular case where $P$ is
a polynomial, the sign change condition can be expressed in terms
of the roots and order of the polynomial $P$, and the
orthogonality properties in terms of moments. For example, for
each $m=0,1,\ldots$ there exists a distributional transformation
which is defined using the Hermite polynomial of order $m$ as the
biasing function, and whose domain are those distributions whose
first $2m$ moments match those of the mean zero normal; the case
$m=1$ corresponds to the zero bias transformation of
\cite{goldrein}.

Theorem \ref{biastheorem} in Section \ref{section2} shows
distributional transformations exist in great generality. In
Section \ref{section3} we find that there is considerable
additional structure for the families of transformations induced
by orthogonal polynomial systems, especially those corresponding
to families of distributions which are closed under addition of
independent variables. Corresponding to the Normal, Gamma,
Poisson, Binomial and Beta-type distributions, in Section
\ref{special-orth-polys} we study the family of transformations
defined using the Hermite, Laguerre, Charlier, Krawtchouk, and
Gegenbauer polynomials, and obtain high order Stein type
characterizing equations.

Our work here is in the spirit of \cite{DZ}, where other
fundamental connections between Stein equations and orthogonal
polynomials were first described. The approach in \cite{DZ} is
iterated in \cite{schoutens}, combining it with well-known
connections between orthogonal polynomials and birth and death
processes, and used as in \cite{holmes} to describe solutions of
Stein equations.

We first review some well known facts regarding the size bias
transformation, which is the simplest and best known of all these
distributional transformations. For non-negative $X$ with $0<EX =
\mu < \infty$, the $X$-size biased distribution $X^s$ is defined
by the characterizing equation \beq \label{size-bias} E XF(X) =
\mu EF(X^s) \quad \mbox{for all $F \in {\cal F}^0(X)$.} \enq One
key feature of the sized bias transformation is the following. If
$X_1,\ldots,X_n$ are independent non-negative variables with
finite positive expectations $EX_i = \mu_i$ and
$$
W = \sum_{i=1}^n X_i,
$$
then a variable with the $W$-size biased distribution can be
constructed by replacing a variable $X_i$, chosen with probability
proportional to $\mu_i$, by an independent variable $X_i^s$ having
the $X_i$-size biased distribution. In other words, letting
$$
P(I=i) = \frac{\mu_i}{\sum_{j=1}^n \mu_j}
$$
be independent of $X_1,\ldots,X_n$, the variable \beq
\label{s-bias} W^s = W - X_I + X_I^s \enq has the $W$-size biased
distribution. Letting $x^+=\max(0,x)$, size biasing is the case of
(\ref{Gbias-intro}) with biasing function $P(x)=x^+$. This
transformation is of order zero, as there are $m=0$ sign changes
of $x^+$ on ${\bf R}$, and has $\alpha = EX^+$; when $X \ge 0$ we
have $X^+=X$ resulting in the usual characterization
(\ref{size-bias}).

The zero bias transformation \cite{goldrein} was motivated by the
similarity between the size bias transformation and the Stein
equation \cite{st72} for the mean zero normal distribution. In
particular, Stein's identity says that $Z \sim {\cal
N}(0,\lambda)$ if and only if
\bea
\label{char-norm} EZF(Z)=\lambda EF'(Z) \quad \mbox{for all $F \in
{\cal F}^1(Z)$}.
\ena
Comparing (\ref{char-norm}) to (\ref{size-bias}), for a mean zero, positive
variance $\lambda$ variable $X$, we say that $X^z$ has the
$X$-zero biased distribution if
\begin{equation}
\label{cheq} EXF(X)=\lambda EF'(X^z) \quad \mbox{for all $F \in
{\cal F}^1(X)$}.
\end{equation}
Note that (\ref{cheq}) for zero biasing is the same as
(\ref{size-bias}) for size biasing, but with variance replacing
mean, and $F'$ replacing $F$. That the normal distribution with
variance $\lambda$ is the unique fixed point of the zero bias
transformation follows immediately from the characterization
(\ref{char-norm}). It was shown in \cite{goldrein} that the zero
bias distribution $X^z$ exists for all $X$ that have mean zero and
finite positive  variance. Its existence follows also from Theorem
\ref{biastheorem}, as the special case of (\ref{Gbias-intro}) for
the function $P(x)=x$, having $m=1$ sign changes on ${\bf R}$, and
$\alpha$ equal to the variance $\lambda$ of $X$.

The zero bias transformation was introduced and used in
\cite{goldrein} to obtain bounds of order $n^{-1}$ in normal
approximations for smooth test functions under third order moment
conditions, in the presence of dependence induced by simple random
sampling. In \cite{hie} it is used to provide bounds to the normal
distribution for hierarchical sequences generated by the iteration
of a so called averaging function, in \cite{zsm} for normal
approximation in combinatorial central limit theorems with random
permutations having distribution constant over cycle type, and in
\cite{goldrein2} the extension of the zero bias transformation to
higher dimension is considered.

The zero bias transformation enjoys a property similar to
(\ref{s-bias}) for size biasing. In particular, it was shown in
\cite{goldrein} that a sum of independent mean zero variables with
finite variances can be zero biased by replacing one variable
chosen with probability proportional to its variance by an
independent variable from that summands zero biased distribution.
Precisely, let $X_1,\ldots,X_n$ be independent mean zero variables
with variance $\lambda_i = EX_i^2 >0$,
$$
W=X_1+ \cdots + X_n,
$$
and $I$ a random index, independent of $X_1,\ldots,X_n$ with
distribution
\bea
\label{pick-var} P(I=i) = \frac{\lambda_i}{\sum_{j=1}^n
\lambda_j}.
\ena
Then \beq \label{z-bias} W^z = W - X_I + X_I^z \enq has the
$W$-zero biased distribution, where $X_i^z$ is a variable
independent of $X_j, j \not = i$ having the $X_i$ zero biased
distribution. This construction is extended to the families of
transformations associated with orthogonal polynomial in Theorem
\ref{big}. In particular, in Section \ref{section3} we see that
for higher order transformations sums of independent variables are
transformed by replacing multiple variables chosen according to
some distribution (e.g. multinomial, multivariate hypergeometric)
with independent variables possessing distributions transformed by
the same family.

In Section \ref{section2} we give the moment and sign change
conditions on $P$ which guarantee the existence of the $X-P$
distribution and provide an explicit construction. In Section
\ref{section3} we treat the special case where $P$ is a member of
a family of orthogonal polynomials. The generalization to higher
order of the `replace one variable' zero and size bias
constructions is based on the identity (\ref{ident:p}) expressing
an orthogonal polynomial of a sum as a sum of like polynomials
with summands having no larger order, and is given in Section
\ref{section3}. In Sections \ref{hermite}, \ref{laggie},
\ref{charlie}, \ref{kraw} and \ref{semi} we treat the Hermite,
Laguerre, Charlier, Krawtchouk and Gegenbauer polynomials,
corresponding to the Normal, Gamma, Poisson, Binomial and
Beta-type distributions respectively. Special instances of the
Beta-type distributions we consider are the uniform ${\cal
U}[-1,1]$, the arcsine, and the semi-circle distribution.

\section{Transformations in General}
\label{section2} We begin our study with the following existence
and uniqueness theorem for the types of distributional
transformations under consideration. We say the measurable
function $P$ on ${\bf R}$ is positive on an interval $I$ if $P(x)
\ge 0$ for all $x \in I$ with strict inequality for at least one
$x$, and similarly for $P$ negative on $I$. We say $P$ has exactly
$m=0,1,\ldots$ sign changes if ${\bf R}$ can be partitioned into
$m+1$ disjoint subintervals with non-empty interior such that $P$
alternates sign on successive intervals. Though the choices for
the endpoints of such intervals may be somewhat arbitrary when
there are intervals where $P$ is zero, we will nevertheless say
that a sign change occurs at the interval boundaries; the
uniqueness guaranteed by Theorem \ref{biastheorem} shows that the
$X-P$ biased distribution constructed in the proof of Theorem
\ref{biastheorem} is the same for all interval boundary choices,
and Example \ref{absolute} gives some additional explanation of
this phenomenon in the context of a particular example. We note
that for existence in general, regarding boundedness, the
orthogonality conditions required by Theorem \ref{biastheorem} are
only relative to $P$ and required only up to a finite order; such
conditions may not impose boundedness on any of the power moments
of $X$, as illustrated in Example \ref{absolute}.

\begin{theorem}
\label{biastheorem} Let $X$ be a random variable, $m \in
\{0,1,2,\ldots\}$ and $P$ a measurable function with exactly $m$
sign changes, positive on its rightmost interval and
\begin{equation}
\label{orth} \frac{1}{m!}EX^k P(X)= \alpha \delta_{k,m} \quad
k=0,\ldots,m,
\end{equation}
with $\alpha >0$. Then there exists a unique distribution for a
random variable $X^{(P)}$ such that
\begin{equation}
\label{Gbias} E P(X) F(X) = \alpha E F^{(m)}(X^{(P)}) \quad
\mbox{for all $F \in {\cal F}^m(P)$.}
\end{equation}
\end{theorem}
Theorem \ref{biastheorem} says that $X$ is in the domain of the
distributional transformation of order $m$ defined using the
`biasing' function $P$ having $m$ sign changes when the powers of
$X$ smaller than $m$ are orthogonal to $P(X)$ in the $L^2(X)$
sense, that is, when $P(X) \in \{1,X,\ldots,X^{m-1} \}^\perp$, and
$EX^m P(X) > 0$. As noted above, the existence of both the size
and zero bias transformations are both special cases.

\noindent {\bf Proof of Theorem \ref{biastheorem}.} We give an
explicit construction of the variate $X^{(P)}$. By replacing $P$
by $P/\alpha$, it suffices to prove the theorem for $\alpha=1$.
Label the points where the $m$ sign changes of $P$ occur as
$r_1,\ldots,r_m$, and let
\bea
\label{Qpoly}
Q(x)=\prod_{i=1}^m (x-r_i),
\ena
adopting the usual convention that an empty product is 1. By
construction $Q$ and $P$ have the same sign, so letting $\mu_X$
denote the distribution of $X$,
\bea
\label{dmuy} d\mu_Y(y)=\frac{1}{m!}Q(y)P(y) d\mu_X(y)
\ena
is therefore a measure, and since (\ref{orth}) with $k=m$ implies
that $EQ(X)P(X)=m!$, a probability measure. Now with $Y$ and
$\{U_i\}_{i \ge 1}$ mutually independent with $Y$ having
distribution $\mu_Y$ and $U_j$ having distribution function $u^i$
on $[0,1]$, with $r_0=Y$ and $r_{m+1}=0$, we claim that
\bea
\label{constructm} X^{(P)} = \sum_{k=1}^{m+1} \left(\prod_{i=k}^m
U_i \right)(r_{k-1}-r_k)
\ena
satisfies (\ref{Gbias}), thus proving the existence of the $X-P$
biased distribution.

We begin by noting that for any $F$ for which either side below
exists,
\begin{equation}
\label{rnxhat} EF(Y)=\frac{1}{m!}EF(X)Q(X)P(X),
\end{equation}
and so for $k=0,\ldots,m$, letting
$$
R_k(x)=\prod_{i=k+1}^m(x-r_i),
$$
a polynomial of degree $m-k$, by (\ref{rnxhat}) and (\ref{orth})
we have
\bea
\label{Rpoly-inverse-expect} E(1/\prod_{i=1}^k(Y-r_i))=
\frac{1}{m!}ER_k(X)P(X) = \delta_{m-k,m}.
\ena

We show the claim by induction. In particular, for $k \ge 1$
letting
\bea
\label{defVW} V_k= \prod_{i=k}^m U_i,  \quad
W_k=\sum_{j=k}^{m+1}V_j(r_{j-1}-r_j),
\ena
and taking $X^{(P)}$ as in (\ref{constructm}), we show that for
all $F \in {\cal C}_c^\infty$, the collection of infinitely
differentiable functions with compact support, and $k=0,\ldots,m$
\begin{equation}
\label{inductp} EF^{(m)}(X^{(P)})=k! \, E \left\{
\frac{F^{(m-k)}(V_{k+1}(Y-r_{k+1})+W_{k+2})}
{V_{k+1}^k\prod_{i=1}^k(Y-r_i)}\right\}.
\end{equation}
We see the expectation on the right hand exists since $F$ and all
its derivatives are bounded, $V_{k+1}$ is independent of $Y$ for
all $k$, $EU_i^{-k}<\infty$ for $i \ge k+1$, and use of
(\ref{Rpoly-inverse-expect}).

The case $k=0$ is the statement that $X^{(P)}=V_1(Y-r_1)+W_2$,
which follows from definitions (\ref{constructm}) and
(\ref{defVW}). Assume (\ref{inductp}) holds for some $0 \le k<m$.
Using $V_{k+1}=U_{k+1}V_{k+2}$ in (\ref{inductp}) and taking
expectation over $U_{k+1}$, with density
 $(k+1)u_{k+1}^k$, we obtain
\beas
\lefteqn{
EF^{(m)}(X^{(P)})}\\
&=&(k+1)!E \int_0 ^1 \left\{
\frac{F^{(m-k)}(u_{k+1}V_{k+2}(Y-r_{k+1})+W_{k+2})} {u_{k+1}^k
V_{k+2}^k\prod_{i=1}^k (Y-r_i)}\right\}u_{k+1}^kdu_{k+1}.
\enas
Cancelling $u_{k+1}^k$ and integrating, we obtain \beq
\label{toobig}
(k+1)!E\left\{\frac{F^{(m-(k+1))}(V_{k+2}(Y-r_{k+1})+W_{k+2})-F^{(m-(k+1))}(W_{k+2})}
{V_{k+2}^{k+1}\prod_{i=1}^{k+1} (Y-r_i)}\right\}. \enq

Using the independence of $V_{k+2}$ and $Y$ for any $k$, and that
$W_{k+2}$ is independent of $Y$ for all $k \ge 0$, the second term
in the expectation (\ref{toobig}) vanishes by
(\ref{Rpoly-inverse-expect}), since $k+1 \ge 1$. The induction is
completed by noting that definitions (\ref{defVW}) give that
$V_{k+2}(Y-r_{k+1})+W_{k+2}=V_{k+2}(Y-r_{k+2})+W_{k+3}$.

Now applying (\ref{inductp}) for $k=m$ and using $V_{m+1}=1,
W_{m+2}=0$ and $r_{m+1}=0$ we obtain
$$
EF^{(m)}(X^{(P)})=m!E\left\{\frac{F(Y)}{Q(Y)}\right\}=EP(X) F(X)
$$
by (\ref{rnxhat}). That is, the equality in (\ref{Gbias}) holds
for all $F \in {\cal C}_c^\infty$.

For $F \in {\cal F}^m(X)$, by replacing $F$ by
$$
F(x)-\sum_{j=0}^{m-1}\frac{F^{(j)}(0)}{j!}x^j
$$
if necessary, we may assume, in light of (\ref{orth}), that
$F^{(j)}(0)=0$ for $j=0,\ldots,m-1$, and hence,
$$
\mbox{with} \quad I f=\int_0^x f, \quad F(x)=I^m f \quad \mbox{for
some measurable function $f$.}
$$
Since $F=F_1-F_2$ where $F_1(x) =I^m f^+$ and $F_2(x) =I^m f^-$,
it suffices by linearity to consider $f \ge 0$. Letting $0 \le f_n
\uparrow f$ we have $I^m f_n = F_n \uparrow F$, and hence the
equality in (\ref{Gbias}) holds for $F \in {\cal F}^m(X)$ using
the monotone and dominated convergence theorems on the right and
left sides of (\ref{Gbias}), respectively.

The distribution $X^{(P)}$ is unique since (\ref{Gbias}) holds for
all $F \in {\cal C}_c^\infty$, which is separating. $\qed$

The existence of the $X^{(P)}$ distribution also follows from the
Riesz representation theorem upon demonstrating the positivity of
the linear operator $T$ defined by
$$
Tf = E  P(X) F(X) \quad \mbox{with $F(x)= I^m f$}
$$
over $f \in {\cal C}_c^0$, the space of continuous functions with
compact support. The signed measure $d\mu=P d\mu_X$ has the
property $\int x^j d\mu = EX^j P(X)=0$ for $j=0,1,\ldots,m-1$, and
now the sign change property of $P$ allows us, when on the finite
interval $[a,b]$, to invoke Theorem 5.4 in Chapter XI of \cite{KS}
(see also Example 1.4 in Chapter XI) to conclude $T$ is positive
and hence $Tf=\int_a^b f d\mu^{(m)}$ for some measure $\mu^{(m)}$,
which is a probability measure since $EX^m P(X)=m!$. This argument
is similar to the one used in \cite{goldrein} to prove the
existence of the zero bias distribution for a mean zero, finite
variance $X$ by noting that when $f \ge 0$ the function $F=If$ is
non-decreasing, and hence $X$ and $F(X)$ are positively
correlated, and so the operator
$$
Tf=EXF(X) \ge EX EF(X) =0
$$
is positive.

\begin{example}
\label{absolute} Consider the application of Theorem
\ref{biastheorem} where $P(x)$ has exactly $m=1$ sign change at
$r_1=0$. Then for the non constant $X$ to be in the domain of the
transformation characterized by
\bea
\label{m=1-trans} EP(X)F(X) = \alpha EF'(X^{(P)})
\ena
we require $EP(X)=0$ and $\alpha=EXP(X)>0$. We have $Q(x)=x$ in
(\ref{Qpoly}) and, recalling the $X$ variable in the proof was
rescaled to have $\alpha=1$, the $Y$ distribution in (\ref{dmuy})
is
$$
d\mu_Y(y)=xP(x)d\mu_X(y)/\alpha.
$$
From (\ref{constructm}) with $m=1$, $r_0=Y, r_2=0$ and $U_j$ with
distribution function $u^j$ on $[0,1]$
$$
X^{(P)} = \sum_{k=1}^{m+1} \left(\prod_{i=k}^m U_i
\right)(r_{k-1}-r_k) = U_1(r_0-r_1)+(r_1-r_2)=U_1Y.
$$
Hence $X^{(P)}$ is absolutely continuous, and one can directly
verify that its density is given by
\bea
\label{fPdensity}
f^{(P)}(x)=\alpha^{-1}E[P(X); X>x].
\ena
When $\int_0^x P(u)du$ is finite for all $x$ and $c=\int
\exp(-\alpha^{-1}\int_0^x P(u)du) dx < \infty$, the transformation
(\ref{m=1-trans}) has a fixed point at the distribution with
density
$$
f(x)=c^{-1}\exp\left(-\frac{1}{\alpha}\int_0^x P(u) du\right);
$$
for instance, when $P(x)=x$, $f$ is the mean zero normal density
with variance $\alpha$.

Taking $P$ to be the sign function
$$
P(x)={\bf 1}(x>0)-{\bf 1}(x<0)
$$
provides an example of a transformation given by a discontinuous
$P$, and shows that generally the orthogonality conditions may not
reduce to restrictions on the moments of $X$, in particular,
(\ref{orth}) for $k=0$ requires $X$ to have median 0. If in
addition $\alpha=E|X|$ is finite, imposed by (\ref{orth}) for
$k=1$, Theorem \ref{biastheorem} gives that $X$ is in domain of
the transformation characterized by (\ref{m=1-trans}). The density
of the transformed variables are, by (\ref{fPdensity}),
\bea
f^{(P)}(x)=\left\{
\begin{array}{cc}
P(X>x)/E|X| & x >0 \\
P(X<x)/E|X| & x<0.
\end{array}
\right.
\ena
For this choice of $P$ the $Y$ distribution in (\ref{dmuy})
becomes
$$
d\mu_Y(y)=|y|d\mu_X(y)/E|X|,
$$
which is the $|X|$ size biased distribution. Hence, the $X-P$
biased distribution is obtain by multiplying $Y \sim \mu_Y$ by an
independent ${\cal U}[0,1]$ variable. The transformation has a
fixed point at the Laplace distribution with density
$$
f(x)=\frac{1}{2\alpha}exp\left(-\frac{1}{\alpha} |x|\right).
$$

Taking $P(x)={\bf 1}(x>1)-{\bf 1}(x<-1)$ gives a transformation
having domain those variables $X$ with $\alpha=E(|X|{\bf 1}(|X|
> 1))<\infty$ and satisfying
\bea
\label{1median} P(X > 1) = P(X < -1).
\ena
Since $P(x)=0$ in the set $[-1,1]$ the sign change can be said to
occur at point in $(-1,1)$ and the polynomial $Q$ in the proof of
Theorem \ref{biastheorem} can be taken to be
$$
Q(x)=x-r_1 \quad \mbox{for any $r_1 \in (-1,1)$.}
$$
As assured by uniqueness, the distribution constructed in the
proof of Theorem \ref{biastheorem} does not depend on choice of
$r_1$; in fact, in this case (\ref{1median}) implies that the
$dmu_Y$ distribution in (\ref{dmuy}) is the same for all $r_1 \in
(-1,1)$.
\end{example}

\section{Transformations using orthogonal polynomials}
\label{section3} We consider a system of polynominals orthogonal
with respect to a non-trivial family of distributions $Z_\lambda
\sim {\cal L}_\lambda$ indexed by a real parameter $\lambda$.

\begin{condition}
\label{poly-condition} For some $m \ge 0$, the polynomials
$\{P^k_\lambda(x)\}_{0\le k \le m}$ are monic, have degree $k$,
are orthogonal with respect to the distributional family
$Z_\lambda \sim {\cal L}_\lambda$, and satisfy
$E[P^k_\lambda(Z_\lambda)]^2>0$.
\end{condition}
Note that since $P_\lambda^k$ is monic and orthogonal it has $k$
distinct roots and is positive as $x \rightarrow \infty$ (e.g.
\cite{orthpolyref}); furthermore, we have
$$
E Z_\lambda^k P_\lambda^k(Z_\lambda) = E
[P_\lambda^k(Z_\lambda)]^2, \quad k=0,\ldots,m.
$$
When studying transformations using an implicit family of
orthogonal polynomials, we index the transformed distribution by
say, $X_\lambda^{(k)}$, that is, by the parameter $\lambda$ and
order $k$ of the polynomial.

Applying Theorem \ref{biastheorem} in this framework, we obtain
the following
\begin{corollary}
\label{cor1} Let Condition \ref{poly-condition} be satisfied with
$E Z_\lambda^{2m} < \infty$, and for $0 \le k \le m$ set
\bea
\label{alpha} \alpha_\lambda^{(k)}  = \frac{1}{k!} EZ_\lambda^k
P^k_\lambda(Z_\lambda).
\ena
Then for all $X \in {\cal M}^k_\lambda$, where
$$
{\cal M}^k_\lambda=\{X:  EX^j=EZ_\lambda^j,
\quad 0 \le j \le 2k \},
$$
there exists a random variable $X^{(k)}_\lambda$ such that for all
$F \in {\cal F}^k(P_\lambda^k)$
\bea
\label{johor-baru} EP^k_\lambda(X)F(X) = \alpha_\lambda^{(k)}
EF^{(k)}(X^{(k)}_\lambda).
\ena
\end{corollary}

\noindent {\bf Proof:} By Condition \ref{poly-condition} and
orthogonality we have for $0 \le j \le k \le m$,
\beas
\frac{1}{k!}EX^jP^k_\lambda(X) = \frac{1}{k!}EZ_\lambda^j
P^k_\lambda(Z_\lambda) = \frac{1}{k!}EP_\lambda^j(Z_\lambda)
P^k_\lambda(Z_\lambda)= \alpha_\lambda^{(k)} \delta_{j,k}
\enas
using $X \in {\cal M}_\lambda^k$. Now invoke Theorem
\ref{biastheorem}. $\qed$

We say the family of distributions $Z_\lambda$ is closed under
independent addition if for independent $Z_{\lambda_i} \sim {\cal
L}_{\lambda_i},i=1,2$ we have $Z_{\lambda_1}+Z_{\lambda_2} \sim
{\cal L}_{\lambda_1+\lambda_2}$. There is special structure when
the transformation function in Theorem \ref{biastheorem} is a
member of an orthogonal polynomial system corresponding to such a
family. In particular, the following Theorem \ref{big} generalizes
(\ref{s-bias}) and (\ref{z-bias}) in showing how a sum of
independent variables can be $P_\lambda^m$ transformed by
replacing a randomly chosen collection in the sum by variables
with distributions transformed using the same orthogonal
polynomial system.

For $n=1,2,\ldots$, consider a multi-index ${\bf
m}=(m_1,\ldots,m_n)$, and with
$\blambda=(\lambda_1,\ldots,\lambda_n)$ and ${\bf
x}=(x_1,\ldots,x_n)$ let
$$
m=|{\bf m}| = \sum_{i=1}^n m_i,  \quad \lambda =
\sum_{i=1}^n\lambda_i
$$
and set
\bea
\label{def:balpha} \alpha^{({\bf m})}_{\blambda}= \prod_{i=1}^n
\alpha^{(m_i)}_{\lambda_i} \quad \mbox{and} \quad P^{\bf
m}_{\blambda}({\bf x})= \prod_{i=1}^n P^{m_i}_{\lambda_i}(x_i).
\ena

\begin{theorem}
\label{big} Let $Z_\lambda, \lambda>0$ be a family of random
variables closed under independent addition with
$EZ_\lambda^{2m}<\infty$, and suppose the associated orthogonal
polynomials $\{P^k_\lambda(x)\}_{0 \le k \le m}$ satisfies
Condition \ref{poly-condition} and, for some weights $c_{\bf m}$,
the identity
\bea
\label{ident:p} P^m_\lambda(w) &=& \sum_{{\bf m}:|{\bf m}|=m}c_{\bf
m} P^{\bf m}_{\blambda}({\bf x}),
\ena
where $P^{\bf m}_{\blambda}({\bf x})$ is given in
(\ref{def:balpha}) and $w=x_1+ \cdots + x_n$. Then
$\alpha^{(m)}_{\lambda}$ and $\alpha^{({\bf m})}_{\blambda}$
defined in (\ref{alpha}) and (\ref{def:balpha}) respectively,
satisfy
\bea
\label{ident:alpha} \alpha^{(m)}_{\lambda} = \sum_{{\bf m}:|{\bf
m}|=m} c_{\bf m} \alpha^{({\bf m})}_{\blambda},
\ena
and we may consider the variable ${\bf I}$, independent of all
other variables, with distribution
\bea
\label{chimes}
P({\bf I}={\bf m})= c_{\bf m} \frac{\alpha^{({\bf
m})}_{\blambda}}{\alpha^{(m)}_{\lambda}}, \quad |{\bf m}|=m.
\ena
Furthermore, for any positive $\lambda_1,\ldots,\lambda_n$ and
independent variables $X_1,\cdots,X_n$ with
\beas
X_i \in {\cal M}_{\lambda_i}^m \quad \mbox{and} \quad W =
\sum_{i=1}^n X_i,
\enas
the variable
\beas
W^{(m)}_\lambda= \sum_{{\bf m}:|{\bf m}|=m}
(X_i)_{\lambda_i}^{(I_i)}
\enas
has the $W-P^m_\lambda$ distribution.
\end{theorem}

\bigskip
\noindent {\bf Proof.} Since $X_i \in {\cal M}_{\lambda_i}^m$, we
have for $0 \le k \le 2m$, and independent $Z_{\lambda_i} \sim
{\cal L}_{\lambda_i}$ and $Z_\lambda \sim {\cal L}_\lambda$, that
\beas
EW^k=E(\sum_{i=1}^n X_i)^k = E(\sum_{i=1}^n Z_{\lambda_i})^k = E
Z_\lambda^k.
\enas
Hence $W \in {\cal M}_\lambda^m$, and the $W^{(m)}$ distribution
exists by Corollary \ref{cor1}. Equality (\ref{ident:alpha})
follows by multiplying (\ref{ident:p}) by $W^m=(\sum_i X_i)^m$
taking expectation, and using independence and orthogonality.

By (\ref{ident:alpha}), for any $F \in {\cal C}_c^\infty$,
\bea
\label{eq:99} \alpha_\lambda^{(m)} EF^{(m)}(W^{(m)}_\lambda) &=& E
\sum_{\bf m} c_{\bf m} \alpha^{({\bf m})}_{\blambda}
F^{(m)}(W^{(m)}_\lambda).
\ena
Using independence and successively applying the identity
\bea
\label{also-works-discrete}
\alpha_{\lambda_i}^{(m_i)}
EF^{(q)}((X_i)_{\lambda_i}^{(m_i)}+y)= EP_{\lambda_i}^{m_i}(X_i)
F^{(q-m_i)}(X_i+y),
\ena
we see that the right hand side of (\ref{eq:99}) is equal to
\bea
\label{little-india} E \sum_{\bf m} c_{\bf m} P^{\bf
m}_{\blambda}({\bf X}) F(W) = E P^m_\lambda(W) F(W),
\ena
by (\ref{ident:p}). Comparing (\ref{eq:99}) to
(\ref{little-india}) we have
$$
\alpha_\lambda^{(m)} EF^{(m)}(W^{(m)}_\lambda) = E P^m_\lambda(W)
F(W),
$$
for all  $F\in {\cal C}_c^\infty$, and hence $W^{(m)}_\lambda$ has
the $W$-$P_\lambda^m$ biased distribution. $\qed$

For the possibly infinite system of monic polynomials
$\{P_\lambda^m(x)\}$ orthogonal with respect to ${\cal
L}_\lambda$, define the generating function
\bea
\label{phi-gen}
\phi_t(x,\lambda)= \sum_{m \ge 0} P_\lambda^m(x)
\frac{t^m}{m!}.
\ena
Though the constants $\alpha_\lambda^{(m)}$ can be found using
$F(x)=x^m$ in (\ref{johor-baru}), squaring (\ref{phi-gen}) and
taking expectation using orthogonality gives the alternative
method
\bea
\label{alpha-gen}
E[\phi_t(Z_\lambda,\lambda)]^2 = \sum_{m \ge 0} \alpha_\lambda^{(m)} \frac{t^{2m}}{m!}.
\ena

Theorem \ref{not-so-big} applies in the special cases considered
in Sections \ref{hermite} through \ref{kraw}.
\begin{theorem}
\label{not-so-big} If the polynomial generating function
$\phi_t(x,\lambda)$ in (\ref{phi-gen}) satisfies
\bea
\label{phi-mult}
\phi_t(w,\lambda) = \prod_{i=1}^n \phi_t(x_i,\lambda_i)
\ena
for $w=x_1+\cdots+x_n$ and $\lambda=\lambda_1+\cdots+\lambda_n$,
then (\ref{ident:p}), and hence (\ref{chimes}), in Theorem
\ref{big} are satisfied respectively by
$$
c_{\bf m}= {m \choose {\bf m}} \quad \mbox{and} \quad P({\bf
I}={\bf m})= {m \choose {\bf m}} \frac{\alpha^{({\bf
m})}_{\blambda}}{\alpha^{(m)}_{\lambda}}, \quad |{\bf m}|=m.
$$
\end{theorem}

\noindent {\bf Proof:} Rewriting (\ref{phi-mult}) ,
\beas
\sum_{m \ge0 } \frac{t^m}{m!} P_\lambda^m(w)
  &=& \prod_{i=1}^n \sum_{m_i \ge 0} P_{\lambda_i}^{m_i}(x_i)\frac{t^{m_i}}{m_i!}\\
&=&
\sum_{m_1,\cdots,m_n}P_{\blambda}^{\bf m}({\bf x})\frac{t^{m_1+\cdots+m_n}}{m_1!\cdots
m_n!}\\
&=& \sum_{m \ge 0}^\infty \frac{t^m}{m!} \sum_{{\bf m}=m}{m
\choose {\bf m}}P_{\blambda}^{\bf m}({\bf x}),
\enas
giving (\ref{ident:p}) with the values claimed. $\qed$

We also note that squaring (\ref{ident:p}) and taking expectation,
using independence and orthogonality, results in
\bea
\label{cmsquared}
\alpha_\lambda^{(m)}=\sum_{|{\bf m}|=m}{m
\choose {\bf m}}^{-1} c_{\bf m}^2 \alpha_{\blambda}^{\bf m},
\ena
so that the conclusion of Theorem \ref{not-so-big} can also be
seen to hold by equating coefficients of (\ref{ident:alpha}) and
(\ref{cmsquared}) when $\alpha_{\blambda}^{\bf m}$ takes on
sufficiently many values.

\bigskip
We end this section with a result about the potential for iterated
biasing.
\begin{theorem}
Let Condition \ref{poly-condition} be
satisfied, and suppose that the the distributional family at
$Z_\lambda$ is closed under transformation with respect to
$P^k_\lambda(x)$, that is, there exists $\mu(\lambda,k)$ such that
$$
(Z_\lambda)_\lambda^{(k)}= Z_{\mu(\lambda,k)}.
$$
Then if $X \in {\cal M}^m_\lambda$ we have $X^{(k)}_\lambda \in
{\cal M}_{\mu(k,\lambda)}^{m-k}$ for $k \le m$. In particular for
non-negative $j$ with $0 \le k+j \le m$, the distribution
$(X^{(k)}_\lambda)_{\mu(\lambda,k)}^{(j)}$ exists.
\end{theorem}
\noindent {\bf Proof.}
 Let $0 \le j \le 2(m-k)$ and $F(x)=x^{k+j}/(k+j)_k$, where $(x)_k=x(x-1)\cdots (x-k+1)$. Then
\beas
&&\alpha^{(k)}_\lambda E(X^{(k)}_\lambda)^j = \alpha^{(k)}_\lambda
EF^{(k)}(X^{(k)}_\lambda)= EP^k_\lambda(X)F(X)\\ &=&
EP^k_\lambda(Z_\lambda)F(Z_\lambda) = \alpha^{(k)}_\lambda
EF^{(k)}((Z_\lambda)_\lambda^{(k)}) = \alpha^{(k)}_\lambda
E(Z_{\mu(\lambda,k)})^j.
\enas
Thus the first $2(m-k)$ moments of $X^{(k)}_\lambda$ match those
of $Z_{\mu(\lambda,k)}$, and the existence of the distribution
$(X^{(k)}_\lambda)_{\mu(\lambda,k)}^{(j)}$ follows from Corollary
\ref{cor1}. $\qed$

\section{Special Orthogonal Polynomial Systems}
\label{special-orth-polys} In Sections \ref{hermite} - \ref{semi}
we specialize to the classic Hermite, Laguerre, Charlier,
Krawtchouk and Gegenbauer orthogonal polynomial systems,
corresponding to the Normal, Gamma, Poisson, Binomial and a Beta
like family, respectively. All these families correspond to a
collection of orthogonal polynomials satisfying Condition
\ref{poly-condition}, and except for the last case, have a
generating function which satisfies (\ref{phi-mult}). The Normal
and Poisson distributions are fixed points of their associated
transformations. In the Gamma, Binomial and Beta-type cases the
transformations map to the same family, but with a shifted
parameter. For further connections between probability
distributions and such polynomial system generating functions, see
\cite{asai} and \cite{asai2}.

\subsection{Hermite Polynomials}
\label{hermite}
For $\sigma^2=\lambda >0$, define the collection
of Hermite polynomials $\{H^m_\lambda(x) \}_{m \ge 0}$ through the
generating function
\bea
\label{herm:gen} e^{xt-\frac{1}{2}\lambda t^2}=\sum_{m=0}^\infty
H^m_\lambda (x) \frac{t^m}{m!},
\ena
or equivalently, the Rodriguez formula
\begin{equation}
\label{herm:rod} H^m_\lambda(x) = (-\lambda)^m
e^{\frac{x^2}{2\lambda}}
\frac{d^m}{dx^m}e^{-\frac{x^2}{2\lambda}}.
\end{equation}
These polynomials are orthogonal with respect to the normal distribution
${\cal N}(0,\lambda)$ with density $(2 \pi
\lambda)^{-1/2}\exp(-x^2/(2 \lambda))$.

For $F \in {\cal C}_c^\infty$ and $Z_\lambda \sim {\cal N}(0,
\lambda)$, applying the Rodriguez formula (\ref{herm:rod}) we have
\bea
\nonumber EH^m_\lambda(Z_\lambda)F(Z_\lambda) &=& \int_{-\infty}
^\infty (-\lambda)^m e^{\frac{x^2}{2\lambda}}
\left(\frac{d^m}{dx^m} e^{-\frac{x^2}{2\lambda}}\right) F(x)
\frac{e^{-\frac{x^2}{2\lambda}}}{\sqrt{\lambda 2\pi}} dx \\
\nonumber &=& \int_{-\infty} ^\infty  (-\lambda)^m
\left(\frac{d^m}{dx^m} e^{-\frac{x^2}{2\lambda}}\right)
F(x) \frac{1}{\sqrt{\lambda 2\pi}} dx \\
\nonumber &=& \lambda^m   \int_{-\infty} ^\infty  F^{(m)}(x)
\frac{e^{-\frac{x^2}{2\lambda}}}{\sqrt{\lambda 2\pi}} dx\\
\label{rod-method-alpha-normal} &=& \lambda^m EF^{(m)}(Z_\lambda).
\ena
Hence,
$$
(Z_\lambda)_\lambda^{(m)} = Z_\lambda,
$$
that is, for each $m=0,1,\ldots,$ the normal $Z_\lambda \sim {\cal
N}(0,\lambda)$ is a fixed point of the $m^{th}$ order
transformation induced by $H_\lambda^m(x)$.

From (\ref{rod-method-alpha-normal}) we see that
$\alpha_\lambda^{(m)} = \lambda^m$, which we could find
alternatively using (\ref{alpha-gen}) and
\beas
E[e^{Z_\lambda t - \frac{1}{2} \lambda t^2}]^2 = e^{\lambda t^2} =
\sum_{m \ge 0} \lambda^m \frac{t^{2m}}{m!}.
\enas
Now since the generating function (\ref{herm:gen}) satisfies the
conditions of Theorem \ref{not-so-big}, the distribution of the
random index $I$ in Theorem \ref{big} is multinomial
$\mbox{Mult}(m,\blambda)$. For zero biasing and $m=1$, this
multinomial distribution reduces to the `pick an index
proportional to variance' as specified in (\ref{pick-var}).

Lastly, we indicate two ways in which the classical Stein equation
can be generalized to the Hermite case. With ${\cal N}h = Eh(Z)$,
the standard normal expectation of $h$, both the equations
\bea \label{h1}
f'(x)H_1^{m-1}(x)-H_1^m(x)f(x) = h(x) - {\cal N}h
\ena
and
\bea \label{h2}
f^{(m)}(x)-H_1^m(x)f(x) = h(x) - {\cal N}h
\ena
reduce to the usual Stein equation when $m=1$ (see \cite{st72},
\cite{st86})
\bea
\label{classical-stein}
f'(x) - xf(x) = h(x) -
{\cal N}h,
\ena
and in particular the expectations on the left hand sides of each
evaluated at a random variable $W$ are zero for all $f \in
C_c^\infty$ if and only if $W$ is standard normal.

\subsection{Laguerre Polynomials}
\label{laggie} For $\lambda>0$, let $\{L^m_\lambda(x) \}_{m \ge
0}$ be the collection of Laguerre polynomials defined by the
generating function
\bea
\label{lagu:gen} (1+t)^{-\lambda}\exp\left\{ \frac{xt}{1+t}
\right\} = \sum_{m=0}^\infty L^m_\lambda (x) \frac{t^m}{m!},
\ena
or equivalently, the Rodriguez formula
\begin{equation}
\label{lagu:rod} L^m_\lambda(x) = (-1)^m x^{-\lambda +1}e^x
\frac{d^m}{dx^m} x^{\lambda +m -1}e^{-x},
\end{equation}
which are orthogonal with respect to the Gamma distribution with
parameter $\lambda$, having density
$x^{\lambda-1}e^{-x}/\Gamma(\lambda), x>0$.

For $F\in {\cal C}_c^\infty$ and $Z_\lambda$ with this density,
applying the Rodriguez formula (\ref{lagu:rod}) yields
\bea
\nonumber EL^m_\lambda(Z_\lambda)F(Z_\lambda) &=& \int_0 ^\infty
(-1)^m x^{-\lambda +1}e^x \left(\frac{d^m}{dx^m} x^{\lambda +m
-1}e^{-x}\right) F(x) \frac{x^{\lambda-1}e^{-x}}{\Gamma(\lambda)} dx \\
\nonumber
&=& \frac{(-1)^m}{\Gamma(\lambda)} \int_0 ^\infty
\left(\frac{d^m}{dx^m} x^{\lambda +m -1}e^{-x}\right)F(x)dx\\
\nonumber
&=& \frac{\Gamma(\lambda+m)}{\Gamma(\lambda)} \int_0
^\infty
F^{(m)}(x)\frac{x^{\lambda +m -1}e^{-x}}{\Gamma(\lambda+m)}dx\\
\label{alpha-by-rod-gamma} &=& (\lambda)^m EF^{(m)}(Z_{\lambda +
m}),
\ena
where $(\lambda)^m$ is the rising factorial,
$$
(\lambda)^m = \lambda(\lambda+1)\cdots(\lambda+m-1)= \frac{\Gamma(\lambda +m)}{\Gamma(\lambda)}.
$$
Hence
$$
(Z_\lambda)_\lambda^{(m)} = Z_{\lambda+m}.
$$

From (\ref{alpha-by-rod-gamma}) we see that $\alpha_\lambda^{(m)}
= (\lambda)^m$, which we could find alternatively using
(\ref{alpha-gen}) and
\beas
E\left[(1+t)^{-\lambda}\exp\left(\frac{Z_\lambda
t}{1+t}\right)\right]^2 = (1-t^2)^{-\lambda} = \sum_{m \ge 0}
(\lambda)^m \frac{t^{2m}}{m!}.
\enas
Since the generating function (\ref{lagu:gen}) satisfies the
conditions of Theorem \ref{not-so-big}, the random index $I$ in
Theorem \ref{big} has distribution
$$
P(I={\bf m}) = {m \choose {\bf m}} \frac{\prod_{i=1}^n
(\lambda_i)^{m_i}}{(\lambda)^m}=\frac{\prod_{i=1}^n
{\lambda_i+m_1-1 \choose m_i} }{{\lambda + m - 1 \choose m}},
$$
which we recognize as the multivariate hypergeometric distribution
with parameters $m$ and $\lambda_1+m_1-1, \ldots,
\lambda_n+m_n-1$, see \cite{JohnsonKotz}, p.301.

Though the Gamma is not a fixed point of the Laguerre
transformations as the normal is for the Hermites, nevertheless
there exist Stein equations for the Gamma paralleling
(\ref{classical-stein}) for the normal which can be used for
studying distributional approximations for the Gamma family; for
details, see \cite{luk}. In particular, we have the Stein
characterization that $X \sim \Gamma(\lambda, 1)$ if and only if
$$
E (X - \lambda ) f(X)= EX f'(X)
$$
for all smoooth functions $f$. Using that $L_\lambda^1(x) = x -
\lambda$, the $X^{(1)}$ order one Laguerre transformation is
characterized by
$$
E (X -
\lambda) f(X) = \lambda E f'(X^{(1)})
$$ for all smooth functions
$f$. Comparing these two equations we see that $X \sim
\Gamma(\lambda, 1)$ if and only if for all smooth functions $f$,
$$ E X f'(X) =  \lambda E f'(X^{(1)});$$
in other words, $X \sim \Gamma(\lambda, 1)$ if and only if
$X^{(1)}$, the first order Laguerre transformation of $X$, equals
its size bias transformation $X^s$.

\subsection{Charlier Polynomials}
\label{charlie} For $\lambda>0$, let $\{C^m_\lambda(x) \}_{m \ge
0}$ be the collection of Charlier polynomials defined by the
generating function
\bea
\label{char:gen}
e^{-\lambda t}\left( 1+t \right)^x=
\sum_{m=0}^\infty C^m_\lambda (x)  \frac{t^m}{m!},
\ena
or, equivalently, with $(x)_k=x(x-1)\cdots (x-k+1)$, the falling
factorial,
\bea
\label{Cm-pre-Rodriguez}
C^m_\lambda (x) = \sum_{k=0}^m {m \choose k}(x)_k
(-\lambda)^{m-k},
\ena
giving a family orthogonal with respect to the Poisson
distribution ${\cal P}(\lambda)$ with mass function
$e^{-\lambda}\lambda^k/k!, k=0,1,\ldots$. From
(\ref{Cm-pre-Rodriguez}) one can derive the Rodriguez formula
\bea
\label{Cm-Rodriguez} C^m_\lambda (x) =
(-1)^m\Gamma(x+1)\lambda^{m-x}\nabla^m\left(
\frac{\lambda^x}{\Gamma(x+1)} \right),
\ena
where $\nabla f(x)=f(x)-f(x-1)$, the backward difference.

Since the transformations in Theorem \ref{biastheorem} defined
using derivatives of test functions yield absolutely continuous
distributions when $m \ge 1$, no discrete distribution will be a
fixed point. However, parallel to (\ref{Gbias}), for an integer
valued random variable $X$ we can define the discrete $X-P$ biased
distribution via
\bea
\label{Gbiasdiscrete}
E P(X) F(X) = \alpha E \Delta^m F (X^{(m)})
\quad \mbox{for all $F \in {\cal F}_\Delta(P)$,}
\ena
where $\Delta f(x)=f(x+1)-f(x)$, and again suppressing dependence
on $X$,
\beas
{\cal F}_\Delta(P)=\{ F:  {\bf R}\rightarrow {\bf R}: E \vert
P(X)F(X)| < \infty \}.
\enas

That for all $m=0,1\ldots$ the Poisson ${\cal
P}(\lambda)$-distribution is a fixed point
$$
(Z_\lambda)_\lambda^{(m)} = Z_\lambda
$$
of the discrete transformation (\ref{Gbiasdiscrete}) with $P$
replaced by $C^m_\lambda$ can be seen as follows. For $Z_\lambda
\sim {\cal P}(\lambda)$, by the Rodriguez formula
(\ref{Cm-Rodriguez}) and
\bea \label{helpful}
\sum_{k=0}^\infty \nabla^m b_k \cdot a_k = (-1)^m
\sum_{k=0}^\infty b_k \Delta^m a_k,
\ena
we have
\bea
\nonumber EC^m_\lambda(Z_\lambda)F(Z_\lambda) &=&
\sum_{k=0}^\infty e^{- \lambda}
\frac{\lambda^k}{k!} C^m_\lambda (k) F(k)\\
\nonumber &=&  \lambda^m (-1)^m \sum_{k=0}^\infty e^{-\lambda}
\nabla^m
\left(\frac{\lambda^k}{k!} \right) \cdot F(k) \\
\nonumber &=&  \lambda^m \sum_{k=0}^\infty
\frac{e^{-\lambda} \lambda^k}{k!} \Delta^m F(k) \\
\label{rod-method-alpha-charlie} &=&\lambda^m E \Delta^m
F(Z_\lambda).
\ena

From (\ref{rod-method-alpha-charlie}) we see that
$\alpha_\lambda^{(m)} = \lambda^m$, which we could find
alternatively using (\ref{alpha-gen}) and
\beas
E[e^{-\lambda t}\left( 1+t \right)^{Z_\lambda}]^2 = e^{\lambda
t^2} = \sum_{m \ge 0} \lambda^m \frac{t^{2m}}{m!}.
\enas

Using the existence of the Charlier biased distributions and that
(\ref{also-works-discrete}) holds with derivative replaced by
difference, it is easy to see that the argument and hence
conclusion of Theorem \ref{big} holds in this discrete case. Now
since the generating function (\ref{char:gen}) satisfies the
conditions of Theorem \ref{not-so-big}, the distribution of the
random index $I$ in Theorem \ref{big} is multinomial
$\mbox{Mult}(m,\blambda)$, as in the normal case.

As  the order one Charlier polynomial is $C_\lambda^1(x) = x -
\lambda$, Stein characterizations of the form (\ref{h1}) or
(\ref{h2}), with Hermite replaced by Charlier and derivatives
replaced by differences, generalize the Stein equation for the
Poisson distribution with parameter $\lambda$ given in \cite{lhy},
and extensively studied for example in \cite{bhj}.

\subsection{Krawtchouk Polynomials}
\label{kraw} With $\lambda = 1,2,\ldots$ and $p \in (0,1)$ fixed,
let $\{K_\lambda^m(x)\}_{0 \le m \le \lambda}$ be the collection
of Krawtchouk polynomials defined by the generating function
\bea
\label{gen:kraw}
\left(1+qt\right)^x(1-pt)^{\lambda-x}=\sum_{m=0}^\lambda
\frac{t^m}{m!} K_\lambda^m(x)
\ena
where $p+q=1$, giving the family of polynomials orthogonal with
respect to the Binomial ${\cal B}(\lambda,p)$ distribution. In
contrast to the previous examples, the Binomial is not infinitely
divisible and has support on a bounded set.

Following the approach set out in \cite{asai2}, the polynomials
can also be given by the Rodriguez formula
\beas
K_\lambda^m(x) = \frac{(-1)^m m! {\lambda \choose
m} p^{m-x} q^x}{{\lambda \choose x} } \nabla^m\left\{ {{\lambda -
m} \choose x} \left( \frac{p}{q}\right)^x \right\},
\enas
and so for $Z_\lambda \sim {\cal B}(\lambda, p), 0 \le m \le
\lambda$ and bounded $F$,
\beas
\lefteqn{E K^m_\lambda(Z_\lambda) F(Z_\lambda) }\\
&=& E F(Z_\lambda)  \frac{(-1)^m m! {\lambda \choose m}
p^{m-Z_\lambda} q^{Z_\lambda}}{{\lambda \choose Z_\lambda} }
\nabla^m\left\{ {{\lambda - m} \choose Z_\lambda} \left( \frac{p}{q}\right)^{Z_\lambda} \right\}\\
&=& \sum_{k=0}^\lambda {\lambda \choose k} p^k q^{\lambda - k}
F(k) \frac{(-1)^m m! {\lambda \choose m} p^{m-k} q^k}{{\lambda
\choose k} }
\nabla^m\left\{ {{\lambda - m} \choose k} \left( \frac{p}{q}\right)^k \right\}\\
&=&  m!{\lambda \choose m} p^m q^\lambda  (-1)^m
\sum_{k=0}^\lambda F(k) \nabla^m\left\{ {{\lambda - m} \choose k}
\left( \frac{p}{q}\right)^k \right\}.
\enas
Using (\ref{helpful}) and letting $(\lambda)_m$ again be the
falling factorial, we write the last expression as
\beas
&& (\lambda)_m p^m q^\lambda   \sum_{k=0}^\lambda
{{\lambda - m} \choose k}  \left( \frac{p}{q}\right)^k \Delta^m F(k)\\
&=& (\lambda)_m (pq)^m  \sum_{k=0}^\lambda
{{\lambda - m} \choose k} p^k q^{\lambda - m - k} \Delta^m F(k)\\
&=& \alpha_\lambda^{(m)} E \Delta^m F(Z_\lambda^{(m)}),
\enas
yielding
$$
\alpha_\lambda^{(m)}= (\lambda)_m (pq)^m \quad \mbox{and} \quad
(Z_\lambda)_\lambda^{(m)}= Z_{\lambda - m}.
$$
Hence, similar to the Gamma family, the Binomial distribution is
not a fixed point of its own transformational family, but the
transformed distribution is a member of the same family. One can
calculate $\alpha_\lambda^{(m)}$ alternatively using
(\ref{alpha-gen}), (\ref{gen:kraw}), and series expansion of
\beas
E (1+  q t)^{2Z_\lambda} (1-pt)^{2 \lambda - 2 Z_\lambda} =(1 +
pqt^2)^\lambda.
\enas

As for Example \ref{charlie}, the conclusion of Theorem \ref{big}
holds, and since the generating function (\ref{gen:kraw})
satisfies the conditions of Theorem \ref{not-so-big} the
distribution of the random index $I$ in Theorem \ref{big} is given
by
\beas
P({\bf I} = {\bf m}) &=& \frac{ {m \choose {\bf m}}  (pq)^{\sum
m_i} \prod_{i=1}^n (\lambda_i)_{m_i}}{ (\lambda)_m(pq)^m}= {m
\choose {\bf m}} \frac{\prod_{i=1}^n
(\lambda_i)_{m_i}}{(\lambda)_m } = \frac{\prod_{i=1}^n {\lambda_i
\choose m_i}}{{\lambda \choose m}},
\enas
which we recognize as the multivariate hypergeometric distribution
with parameters $m$ and $\lambda_1, \ldots, \lambda_n$, see
\cite{JohnsonKotz}, p.301.

From \cite{ehm} we have the Stein characterization that $X \sim
{\cal B}(\lambda, p)$ if and only if
$$
p E(\lambda - X) f(X+1) = q E X f(X),
$$
for all functions $f$ for which these expectations exist. Using
the first Krawtchouk polynomial is $K_\lambda^1(x) =qx-p(\lambda -
x)$, we obtain that the first order Krawtchouk transformation is
characterized by
$$
q E X f(X) - pE (\lambda - X) f(X) = \lambda p q E \Delta
f(X^{(1)}). $$ Combining these equations yields that $X \sim {\cal
B}(\lambda, p)$ if and only if
\beas
p E(\lambda - X) \Delta f(X) &=& p E(\lambda - X) (f(X+1) - f(X)
)\\
&=& q EXf(X) - p E (\lambda - X)f(X) \\
&=&\lambda p q E \Delta f(X^{(1)}).
\enas
Putting $g(x) = \Delta f(\lambda - x)$ we see that $X \sim {\cal
B}(\lambda, p)$ if and only if  $\lambda - X^{(1)}$ has the
$(\lambda - X)$-size biased distribution, that is, if and only
if
$$
\lambda-X^{(1)} \sim {\cal B}(\lambda-1,q)+1, \quad \mbox{which is
equivalent to} \quad X^{(1)} \sim {\cal B}(\lambda-1,p).
$$

\subsection{Gegenbauer Polynomials}
\label{semi} In this last section, we consider a polynomial system
orthogonal with respect to a continuous distribution with compact
support. For $\lambda > -\frac{1}{2}$,
let (see \cite{asai2})
$$
\alpha_\lambda^{(m)} = \frac{\Gamma(\lambda) \Gamma(2 \lambda + m)
\Gamma(\lambda + 1)}{ 2^{2m} \Gamma(\lambda + m + 1 )
\Gamma(\lambda + m ) \Gamma (2 \lambda) },
$$
and the collection of Gegenbauer polynomials $G^m_\lambda(x)$ be
defined via the Rodriguez formula
$$
G^m_\lambda(x) = (-1)^m  \alpha_\lambda^{(m)} \frac{\Gamma(\lambda
+ \frac{1}{2}) \Gamma(\lambda + m + 1) }{ \Gamma( \lambda+1)
\Gamma(\lambda + m + \frac{1}{2})} (1-x^2)^{\frac{1}{2} - \lambda}
\frac{d^m}{dx^m} \{(1-x^2)^{\lambda + m - \frac{1}{2}}\} .
$$
Then $G^m_\lambda(x)$ are monic, have degree $m$, and satisfy the
orthogonality relation
$$
\frac{1}{m!}EG_\lambda^k(Z_\lambda) G_\lambda^m(Z_\lambda) =
\alpha_\lambda^{(m)} \delta_{k,m} \quad k=0,\ldots,m,
$$
where $Z_\lambda \sim g_\lambda$ with
$$
g_\lambda(x)  = \frac{1}{\sqrt{\pi}} \frac{\Gamma(\lambda +
1)}{\Gamma(\lambda + \frac{1}{2})} (1-x^2)^{\lambda -
\frac{1}{2}}, \quad \vert x \vert \le 1.
$$
In particular Corollary \ref{cor1} obtains, proving the existence
of the family of Gegenbauer transformations. This family of
distribution is a special case of the centered Pearson Type
I-distributions, sometimes also called Beta Type I-distributions,
see \cite{kendallstuart}, p.150. We note that for $\lambda=1/2$ we
obtain the uniform distribution ${\cal U}[-1,1]$, for $\lambda =0$
the arcsine law, and for $\lambda =1$ the semi-circle law
\cite{Wigner1}, \cite{Wigner2}.

Considering the action of the $G_\lambda^m$ transformation on
$Z_\lambda \sim g_\lambda$, for $F \in {\cal C}_c^\infty$ we have
\beas
\lefteqn{E G^m_\lambda(Z_\lambda) F(Z_\lambda)}\\
&=& \frac{1}{\sqrt{\pi}} \frac{\Gamma(\lambda + 1)}{\Gamma(\lambda
+ \frac{1}{2})} (-1)^m \alpha_\lambda^{(m)}  \frac{\Gamma(\lambda
+ \frac{1}{2}) \Gamma(\lambda+m+1 ) }{ \Gamma( \lambda+1)
\Gamma(\lambda + m+ \frac{1}{2})} \int_{-1}^1 F(x)
\frac{d^m}{dx^m} \{(1-x^2)^{\lambda + m -
\frac{1}{2}}\} \\
&=& \frac{1}{\sqrt{\pi}} \alpha_\lambda^{(m)}  \frac{
\Gamma(\lambda+m+1 ) }{ \Gamma(\lambda + m+ \frac{1}{2})}
\int_{-1}^1 F^{(m)}(x) (1-x^2)^{\lambda + m -
\frac{1}{2} } \\
&=& \alpha_\lambda^{(m)}  E F^{(m)}(Z_{\lambda+m}),
\enas
yielding
$$
(Z_\lambda)_\lambda^{(m)}=Z_{\lambda+m}.
$$
Thus, for $\lambda =0$ we obtain that the first order
Gegenbauer transformation of the arcsine distribution is the
semi-circle law.

Lastly we note that since the above Beta-type distributions are
not invariant under addition, Theorem \ref{big} and its
construction do not apply. However, as $G_\lambda^1(x) = x$, we
recognize the first order Gegenbauer transformation as the
zero-bias transformation, so that for sums of independent random
variables the construction given in (\ref{z-bias}) applies.

\bigskip
{\bf Acknowledgement.} The authors would like to thank the
organizers of the program `Stein's method and applications: A
program in honor of Charles Stein' and the Institute of
Mathematical Sciences in Singapore for their most generous
hospitality, and the excellent meeting where this work was
completed. Also we would like to thank an anonymous referee for very helpful comments.

\end{document}